\documentclass{article}
\usepackage{array,amsmath}
\numberwithin{equation}{section}
\usepackage{indentfirst}

\begin{document}

\newtheorem{thm}{Theorem}[section]
\newtheorem{cor}[thm]{Corollary}
\newtheorem{prop}[thm]{Proposition}
\newtheorem{conj}[thm]{Conjecture}
\newtheorem{lem}[thm]{Lemma}
\newtheorem{Def}[thm]{Definition}
\newtheorem{rem}[thm]{Remark}
\newtheorem{prob}[thm]{Problem}
\newtheorem{ex}{Example}[section]

\newcommand{\be}{\begin{equation}}
\newcommand{\ee}{\end{equation}}
\newcommand{\ben}{\begin{enumerate}}
\newcommand{\een}{\end{enumerate}}
\newcommand{\beq}{\begin{eqnarray}}
\newcommand{\eeq}{\end{eqnarray}}
\newcommand{\beqn}{\begin{eqnarray*}}
\newcommand{\eeqn}{\end{eqnarray*}}
\newcommand{\bei}{\begin{itemize}}
\newcommand{\eei}{\end{itemize}}

\newcommand{\pa}{{\partial}}
\newcommand{\V}{{\rm V}}
\newcommand{\R}{{\bf R}}
\newcommand{\K}{{\rm K}}
\newcommand{\e}{{\epsilon}}
\newcommand{\tomega}{\tilde{\omega}}
\newcommand{\tOmega}{\tilde{Omega}}
\newcommand{\tR}{\tilde{R}}
\newcommand{\tB}{\tilde{B}}
\newcommand{\tGamma}{\tilde{\Gamma}}
\newcommand{\fa}{f_{\alpha}}
\newcommand{\fb}{f_{\beta}}
\newcommand{\faa}{f_{\alpha\alpha}}
\newcommand{\faaa}{f_{\alpha\alpha\alpha}}
\newcommand{\fab}{f_{\alpha\beta}}
\newcommand{\fabb}{f_{\alpha\beta\beta}}
\newcommand{\fbb}{f_{\beta\beta}}
\newcommand{\fbbb}{f_{\beta\beta\beta}}
\newcommand{\faab}{f_{\alpha\alpha\beta}}

\newcommand{\pxi}{ {\pa \over \pa x^i}}
\newcommand{\pxj}{ {\pa \over \pa x^j}}
\newcommand{\pxk}{ {\pa \over \pa x^k}}
\newcommand{\pyi}{ {\pa \over \pa y^i}}
\newcommand{\pyj}{ {\pa \over \pa y^j}}
\newcommand{\pyk}{ {\pa \over \pa y^k}}
\newcommand{\dxi}{{\delta \over \delta x^i}}
\newcommand{\dxj}{{\delta \over \delta x^j}}
\newcommand{\dxk}{{\delta \over \delta x^k}}

\newcommand{\px}{{\pa \over \pa x}}
\newcommand{\py}{{\pa \over \pa y}}
\newcommand{\pt}{{\pa \over \pa t}}
\newcommand{\ps}{{\pa \over \pa s}}
\newcommand{\pvi}{{\pa \over \pa v^i}}
\newcommand{\ty}{\tilde{y}}
\newcommand{\bGamma}{\bar{\Gamma}}

\font\BBb=msbm10 at 12pt
\newcommand{\Bbb}[1]{\mbox{\BBb #1}}

\newcommand{\qed}{\hspace*{\fill}Q.E.D.}  

\title{The characterizations on a class of weakly weighted Einstein-Finsler metrics\footnote{The first author is partially supported by the National Natural Science Foundation of China (11871126, 12141101) and the Science Foundation of Chongqing Normal University (17XLB022)}}
\author{ Xinyue Cheng\footnote{Corresponding author}, Hong Cheng  and Pengsheng Wu} 
\date{}

\maketitle

\begin{abstract}
In this paper, we study the weakly weighted Einstein-Finsler metrics. First, we show that weakly weighted Einstein-Kropina metrics must be of isotropic S-curvature with respect to the Busemann-Hausdorff volume form under a certain condition about the weight constants. Then we characterize weakly weighted Einstein-Kropina metrics completely via their navigation expressions or via $\alpha$ and $\beta$ respectively.
{\bf MR(2020) Subject Classification:}  53B40, 53C60
\end{abstract}

\section{Introduction}

In Riemannian geometry, the $\infty$-Bakry-Emery Ricci Tensor on a smooth Riemannian metric space $\left(M, g, e^{-f}d{\rm vol}_{g}\right)$ is defined as
\[
\operatorname{Ric}_{\infty}=\operatorname{Ric}+\operatorname{Hess}f.
\]
Here, $M$ is a complete $n$-dimensional Riemannian manifold with metric $g$, and $f$ is a smooth real valued function on $M$,  $d\mathrm{vol}_{g}$ is the standard  Riemannian volume form of $g$ on $M$. The equation $\operatorname{Ric}_{\infty}=\lambda g$ for some constant $\lambda$ is exactly the gradient Ricci soliton equation, which plays an important role in the theory of Ricci flow. More generally, for $N \in {R} \backslash\{n\}$, the $N$-weighted Ricci curvature tensor associated with the measure $dV = e^{-f}d{\rm vol}_{g}$ is defined as
\[
{\rm Ric}_{N}= {\rm Ric}+ {\rm Hess}f-\frac{df \otimes df}{N-n}.
\]
The equation $\operatorname{Ric}_{N}=\lambda g$ for some constant $\lambda$ is a special case of the so-called generalized quasi-Einstein equation (\cite{Cat}).

Recently, the study of the weighted Ricci curvature in Finsler geometry has attracted a lot of attentions. Roughly speaking, the weighted Ricci curvatures in Finsler geometry are various combinations of Ricci curvature and S-curvature,  where the S-curvature was first introduced by Z. Shen when he studied Bishop-Gromov volume comparison in Finsler geometry (\cite{Sh}). Let $(M, F, dV)$ be an $n$-dimensional Finsler measure manifold with  $dV = \sigma (x) dx^{1} \cdots dx^{n}$. Let $Y$ be a $C^{\infty}$ geodesic field on an open subset $U \subset M$ and $\hat{g}=g_{Y}.$  Let
$$
dV:=e^{- \psi} d{\rm vol}_{\hat{g}}, \ \ \ d{\rm vol}_{\hat{g}}= \sqrt{{det}\left(g_{i j}\left(x, Y_{x}\right)\right)}dx^{1} \cdots dx^{n}.
$$
It is easy to see that $\psi$ is given by
$$
\psi (x)= \ln \frac{\sqrt{\operatorname{det}\left(g_{i j}\left(x, Y_{x}\right)\right)}}{\sigma(x)}=\tau\left(x, Y_{x}\right),
$$
which is just the distortion along $Y_{x}$ at $x\in M$ (\cite{ChernShen}).
Let $y := Y_{x}\in T_{x}M$ (that is, $Y$ is a geodesic extension of $y\in T_{x}M$). Then, by the definitions of the S-curvature and the Hessian (\cite{shen1}\cite{shen2}), we have
\beqn
&& {\bf S}(x, y)=y[\tau(x, Y_{x})]= d \psi (y),\\
&& \dot{\bf S}(x, y)=y[{\bf S}(x, Y)]=y[Y(\psi)]= {\rm Hess} \psi (y).
\eeqn
Then, for $N \in R \backslash\{n\},$  we define the N-weighted Ricci curvature (\cite{ChSh})
\be
{\rm Ric}_{N}(y)= {\rm Ric}(y)+ {\rm Hess}  \psi(y)-\frac{d \psi(y)^{2}}{N-n}.   \label{weRicci3}
\ee
Obviously, (\ref{weRicci3}) is an analogue of the N-weighted Ricci curvature in Riemannian geometry. As the limits of $N \rightarrow \infty$ ,  we define the $\infty$-weighted Ricci curvature as follows.
\be
{\rm Ric}_{\infty}(y): = {\rm Ric}(y)+ {\rm Hess} \psi(y). \label{inftyweRic}
\ee
Actually, (\ref{weRicci3}) and (\ref{inftyweRic}) are just the weighted Ricci curvatures introduced by Ohta in \cite{Oh1} as follows
\beq
&& {\rm Ric}_{N}(y) = {\rm Ric}(y) +\mathbf{\dot{S}}-\frac{1}{N-n}\mathbf{S}^{2}, \label{weNRic}\\
&& {\rm Ric}_{\infty}(y) = {\rm Ric}(y) +\mathbf{\dot{S}}. \label{weinfRic}
\eeq

There is another weighted Ricci curvature in Finsler geometry,  namely, the projective Ricci curvature
\be
{\rm PRic}(y) = {\rm Ric}(y) +(n-1)\left[\frac{\mathbf{\dot{S}}}{n+1}+\frac{\mathbf{S}^{2}}{(n+1)^{2}}\right].
\ee
The projective Ricci curvature is a projective invariant when the volume form $dV$ is fixed (\cite{CSM}\cite{SS}).

Recently, Z. Shen and R. Zhao define the notion of the $(a,b)$-weighted Ricci curvature in Finsler geometry (\cite{SZ}), which we will call the generalized weighted Ricci curvature in this paper. Concretely,  a generalized weighted Ricci curvature with weight constants $a$ and $c$ is defined as follow
\be
{\rm Ric}_{a,c} (y) = {\rm Ric}(y) + a {\dot{\bf S}}- c {\bf S}^{2}. \label{gwRic}
\ee
Obviously, the weighted Ricci curvatures  ${\rm Ric}_{N}$, ${\rm Ric}_{\infty}$ and the projective Ricci curvature ${\rm PRic}$ are all special generalized weighted Ricci curvatures with special weight constants.  The generalized weighted Ricci curvature with weight constants $a$ and $c$ can also be written in the following form
\be
{\rm Ric}_{a,c}= {\rm PRic} -\frac{\kappa}{n+1}(\mathbf{\dot{S}}+\frac{4}{n+1}\mathbf{S}^{2})+\frac{\nu}{(n+1)^{2}}\mathbf{S}^{2}, \label{gwRic2}
\ee
where $\kappa:=(n-1)-a(n+1)$ and $\nu:=3(n-1)-4a(n+1)-c(n+1)^{2}$.

A Finsler metric $F$ on an $n$-dimensional manifold $M$ with a volume form $dV=e^{-(n+1)f}dV_{BH}$ is said to be a weakly weighted Einstein-Finsler metric with weight constants $a$ and $c$
if the generalized weighted Ricci curvature satisfies
\be
{\rm Ric}_{a,c}=(n-1)\left(\frac{3\theta}{F}+\sigma\right)F^{2}, \label{weweEF}
\ee
where $dV_{BH}$ denotes the Busemann-Hausdorff volume form and $\sigma$ is a scalar function and $\theta=\theta_{i}y^{i}$ is a 1-form on $M$. If $F$ satisfies (\ref{weweEF}) with $\theta =0$, then $F$ is called a weighted Einstein-Finsler metric with weight constants $a$ and $c$. In particular, when $a=c=0$, Finsler metrics satisfying (\ref{weweEF}) are called the weak Einstein metrics (\cite{ChSh0}). When $a=1$, $c=0$ and $\theta =0$, (\ref{weweEF}) becomes $\operatorname{Ric}_{\infty}= (n-1)\sigma F^{2}$ and $(M, F, dV)$ is called a Finsler gradient Ricci almost soliton. If $(M, F, dV)$ satisfies $\operatorname{Ric}_{\infty}= (n-1)\sigma F^{2}$ for a constant $\sigma$, it is called a Finsler gradient Ricci soliton. In \cite{MZZ}, Mo-Zhu-Zhu find the sufficient and necessary conditions for a Randers measure space to be a Finsler gradient Ricci soliton. More generally, Shen-Zhao find that weakly weighted Einstein-Randers metrics must be of isotropic S-curvature with respect to the Busemann-Hausdorff volume form when $\nu \neq 0$. Further, they classify this class of Randers metrics completely via their navigation expressions or $\alpha$ and $\beta$ respectively (\cite{SZ}).

In this paper, we will specialize in studying weakly weighted Einstein-Kropina metrics. We first find that every weakly weighted Einstein-Kropina metric must be of isotropic S-curvature with respect to the Busemann-Hausdorff volume form if $\nu \neq 0$ (see Theorem \ref{mainth1}). Then we characterize weakly weighted Einstein-Kropina metrics completely via their navigation expressions and via $\alpha$ and $\beta$  respectively when $\nu\neq 0$ (see Theorem \ref{mainth3} and Theorem \ref{mainth2}). In particular, when  $\nu \neq 0$ and $S$-curvature with respect to the Busemann-Hausdorff volume form is isotropic, we prove that a Kropina metric $F$ with navigation data $(h, W)$ is a weakly weighted Einstein metric if and only if $h$ is a  weighted Einstein-Riemann metric satisfying (\ref{formula13}) with respect to volume form $dV_{f}=e^{-(n+1)f}d\mathrm{vol}_{h}$. We also characterize weakly weighted Einstein-Kropina metrics when $\nu =0$ and $\kappa \neq 0$ (see Theorem \ref{mainth4}) or when $\nu =0$ and $\kappa = 0$ (see Theorem \ref{mainth5}) respectively.

\section{Preliminaries}\label{Sec2}

Let $F$ be a Finsler metric on an $n$-dimensional manifold $M$ and $G^{i}$ be the geodesic coefficients of $F$,
which are defined by
\begin{equation}\label{geodesic1}
 G^{i}:=\frac{1}{4}g^{il}\left\{[F^{2}]_{x^{k}y^{l}}y^{k}-[F^{2}]_{x^{l}}\right\}.
\end{equation}
For any $x\in M$ and $y\in T_{x}M\backslash \{0\}$, the Riemann curvature $\mathbf{{R}}_{y}:=R^{i}_{ \ k}\frac{\partial}{\partial x^{i}}\otimes dx^{k}$ is defined by
\begin{equation}
R^{i}_{\ k}:=2\frac{\partial G^{i}}{\partial x^{k}}-\frac{\partial^{2}G^{i}}{\partial x^{m}\partial y^{k}}y^{m} +2G^{m}\frac{\partial^{2}G^{i}}{\partial y^{m}\partial y^{k}}-\frac{\partial G^{i}}{\partial y^{m}}\frac{\partial G^{m}}{\partial y^{k}}. \label{Riemann1}
\end{equation}
The Ricci curvature is defined as the trace of the Riemann curvature, that is,
\begin{equation}
{\rm Ric}(y):=R^{m}_{\ m}. \label{Ricci1}
\end{equation}

Let $(M, F, dV)$ be a Finsler measure manifold with a volume form $dV=\sigma(x) dx^{1} \cdots  dx^{n}$. For each $y \in T_{x}M\setminus \{0\}$,
the quantity
\be
\tau(x,y):= \ln \frac{\sqrt{det\left(g_{ij}(x,y)\right)}}{\sigma (x)} \label{distortion}
\ee
is called the distortion of $F$. Further, let $c = c(t)$ be the geodesic with $c(0)=x$ and $\dot{c}(0)=y$. The S-curvature ${\bf S}$ and its change rate $\dot{\bf S}$ along geodesic $c$ are defined by
\be
\mathbf{S}(x,y):=\left.\frac{d}{dt}\left[\tau(c(t),\dot{c}(t))\right]\right|_{t=0}, \ \ \dot{\bf S}(x,y):=\left.\frac{d}{dt}\left[\mathbf{S}(c(t),\dot{c}(t))\right]\right|_{t=0} \label{SandSdot}
\ee
respectively. In short, we have
\be
\mathbf{S}= \tau_{|m}(x, y) y^{m},\ \ \ \dot{\bf S}={\bf S}_{|m}(x,y)y^{m}, \label{SandSdot1}
\ee
where $`` |"$ denotes the horizontal covariant derivative with respect to $F$.

The Busemann-Hausdorff volume form $d V_{BH}=\sigma_{BH}(x) d x^{1} \cdots d x^{n}$ on a Finsler manifold $(M, F)$ is defined by (\cite{ChernShen})
\[
\sigma_{B H}(x):= \frac{\operatorname{Vol}\left(\mathbf{B}^{n}(1)\right)}{\operatorname{Vol}\left\{\left(y^{i}\right) \in \mathbf{R}^{n} \mid F(x, y)<1\right\}}.
\]
In the following, we always use ${\bf S}_{BH}$ to denote the S-curvature determined by Busemann-Hausdorff volume form $d V_{BH}$.

Let $f$ be a $C^2$ function on $M$. The Hessian of $f$ can be defined as a map ${\rm Hess}_{F} f:  TM\rightarrow R$ by
\be
{\rm Hess}_{F}f(y):=\frac{d^2}{ds^2}\left(f\circ c\right)|_{s=0}, \ \ \ y\in T_{x}M,\label{hessian1}
\ee
where $c: (-\varepsilon , \varepsilon) \rightarrow M$ is the geodesic with $c(0)=x, \ \dot{c}(0)=y\in T_{x}M$ (see (\cite{shen2})).  In local coordinates,
\beq
{\rm Hess}_{F} f(y)&=& \frac{\pa ^{2}f}{\pa x^{i}\pa x^{j}}(x)\dot{c}^{i}(0)\dot{c}^{j}(0)+\frac{\pa f}{\pa x^{i}}(x)\ddot{c}^{i}(0) \nonumber\\
&=& \frac{\pa ^{2}f}{\pa x^{i}\pa x^{j}}(x)y^{i}y^{j}-2\frac{\pa f}{\pa x^{i}}(x)G^{i}(x,y) \nonumber \\
&=& \left(\frac{\pa ^{2}f}{\pa x^{i}\pa x^{j}}(x)-\frac{\pa f}{\pa x^{m}}\Gamma ^{m}_{ij}(x,y)\right)y^{i}y^{j}. \label{Hessian}
\eeq
Here, $\Gamma ^{k}_{ij}(x,y)$ denote the Chern connection coefficients of $F$, which depends on the tangent vector $y\in T_{x}M$ usually.

Randers metrics and Kropina metrics form an important class of Finsler metrics which are both called $C$-reducible Finsler metrics when $n\geq 3$. Randers metrics and Kropina metrics can be also both expressed as the solution of the Zermelo navigation problem on some Riemannian manifold $(M, h)$ with a vector field $W$ (\cite{CQX}\cite{ChSh0}\cite{RS}). Randers metrics are one of the simplest non-Riemannian Finsler metrics with the form $F=\alpha+\beta$, where $\alpha=\sqrt{a_{ij}(x)y^{i}y^{j}}$ is a Riemannian metric and $\beta = b_{i}(x)y^{i}$ is a 1-form with $\|\beta _{x}\|_{\alpha}<1$ on the manifold.  Kropina metrics are those Finsler metrics in the form $F=\frac{\alpha ^2}{\beta}$. Kropina metrics were first introduced by L. Berwald when he studied the two-dimensional Finsler spaces with rectilinear extremal and were investigated by Kropina (see \cite{VK1}\cite{VK2}). Kropina metrics have important and interesting applications in the theory of thermodynamics. Besides, both of Randers metrics and Kropina metrics play an interesting role in the Krivan problem in ecology (\cite{AIM}).  However,  Randers metrics are regular Finsler metrics, but Kropina metrics are the Finsler metrics with singularity. In fact, Kropina metrics are not classical Finsler metrics, but conic Finsler metrics  defined on the conic domain (\cite{CLY}\cite{Ma}\cite{XQ}\cite{RS})
\[
A =\big\{(x,y) \in TM \ | \ \beta= b_{i}(x)y^{i} >0\big\}\subset TM.
\]

Let
\[
r_{ij}:=\frac{1}{2}(b_{i;j}+b_{j;i}),\ \ \  s_{ij}:=\frac{1}{2}(b_{i;j}-b_{j;i}),
\]
where `` ;" denotes the covariant derivative with respect to the Levi-Civita connection of $\alpha$. Put
\[
r^{i}_{\ j}:=a^{ik}r_{kj},\ \ r_{j}:=b^{i}r_{ij},\ \ r^{i}:=a^{ij}r_{j}, \ \ r:=r_{ij}b^{i}b^{j},
\]
\[
s^{i}_{\ j}:=a^{ik}s_{kj},\ \ s_{j}:=b^{i}s_{ij},\ \ s^{i}:=a^{ij}s_{j}, \ \ e_{ij}:=r_{ij}+b_{i}s_{j}+b_{j}s_{i},
\]
where $(a^{ij}):=(a_{ij})^{-1}$ and $b^{i}:=a^{ij}b_{j}$. Besides, we will use the following notations: $r_{i0}=r_{ij}y^{j},\  r_{00}=r_{ij}y^{i}y^{j}, \ s_{i0}=s_{ij}y^{j}, etc.$. Obviously, $ s_{00}=0.$

Let $F=\frac{\alpha ^2}{\beta}$ be a Kropina metric and $G^{i}$ and $G^{i}_{\alpha}$ be the geodesic coefficients of $F$ and $\alpha$ , respectively. Then we have the following lemmas.

\begin{lem}{\rm (\cite{ZS})} For the Kropina metric $F=\frac{\alpha^{2}}{\beta}$, its geodesic coefficients $G^{i}$ are connected with the geodesic coefficients $G^{i}_{\alpha}$ of $\alpha$  by
\be
G^{i}=G^{i}_{\alpha}+T^{i}, \label{geodesic2}
\ee
where
\be
T^{i}=-\frac{\alpha^{2}}{2\beta}{s^{i}}_{0}+\frac{1}{2b^{2}}\left(\frac{\alpha^{2}}{\beta}s_{0}+r_{00}\right)b^{i}-\frac{1}{b^{2}}\left(s_{0}+\frac{\beta}{\alpha^{2}}r_{00}\right)y^{i}. \label{geodesicT}
\ee
\end{lem}

\begin{lem} {\rm (\cite{XQ}\cite{ZS})} For the Kropina metric $F=\frac{\alpha^{2}}{\beta}$, the Ricci curvature of $F$ is given by
\be
{\rm Ric} = {\rm Ric}^{\alpha}+ T, \label{Ricci2}
\ee
where ${\rm Ric}^{\alpha}$ denotes the Ricci curvature of $\alpha$, and
\beq
T&=&\frac{3(n-1)}{b^{4}F^{2}}r^{2}_{00}+\frac{n-1}{Fb^{4}}\left(2r_{00}s_{0}-4r_{00}r_{0}-4Fr_{0}s_{0}-Fs^{2}_{0}\right) \nonumber\\
&&+\frac{n-1}{b^{2}F}\left(r_{00;0}+Fs_{0;0}+F^{2}s_{k}{s^{k}}_{0}\right)+\frac{1}{b^{4}}\left[(r_{0}+s_{0})^{2}-r(r_{00}+Fs_{0})\right] \nonumber\\
&&+\frac{1}{b^{2}}\Big\{(Fs_{0;k}+r_{00;k})b^{k}-(r_{0;0}+s_{0;0})+(r_{00}+Fs_{0})r^{k}_{\ k}+2nr_{0k}s^{k}_{\ 0} \nonumber \\
&& -Fr_{k}s^{k}_{\ 0}-Fr_{0k}s^{k}-\frac{F^{2}}{2}s^{k}s_{k}\Big\}-Fs^{k}_{\ 0;k}-\frac{F^{2}}{4}s^{j}_{\ k}s^{k}_{\ j}. \label{RicciT}
\eeq
\end{lem}

Kropina metrics can be expressed as the solution of the Zermelo navigation problem on some Riemannian manifold $(M, h)$ with a vector field $W$. Concretely, assume that $h=\sqrt{h_{ij}(x)y^iy^j}$ and  $W=W^{i}\frac{\partial}{\partial x^{i}}$ with $\|W\|_h =1$. Then the metric $F$ obtained by solving the following problem
\be
h\left(x,\frac{y}{F(x,y)} -W_{x}\right) =1     \label{eqa2}
\ee
is a Kropia metric given  by
\be
F=\frac{h^{2}}{2W_{0}}, \label{hwKropina}
\ee
where $W_{0}:=W_{i}y^{i}=  h (y,W_{x})$, $W_{i}:=h_{ij}W^{j}$. In this case, the pair $(h, W)$ is called the navigation data of conic Kropina metric $F=\frac{\alpha ^2}{\beta}$ and
\be
\alpha =\frac{b}{2}h, \ \ \ \beta=\frac{b^2}{2}W_0.\label{eqa4}
\ee
Here, $b :=\|\beta\|_\alpha$ denotes the norm of $\beta$ with respect to $\alpha$.  In fact, $F$ given by (\ref{hwKropina}) is a conic Kropina metric defined on the conic domain
\[
A =\big\{(x,y) \in TM \ | \ h (y,W_{x}) >0\big\}= \big\{(x,y) \in TM \ | \ \beta= b_{i}(x)y^{i} >0\big\} \subset TM.
\]
 Conversely, given a conic Kropina metric $F=\frac{\alpha^2}{\beta}$, put
\be
h_{ij}=\frac{4}{b^2}a_{ij},  \ \ \ W^{i}=\frac{1}{2}b^{i}. \label{NaKro}
\ee
Then we can get a Riemannian metric $h$ and a vector field $W$ with $\|W\|_{h}=1$ from (\ref{NaKro}) and $F$ is just given by (\ref{eqa2}) for $h$ and $W$.
Thus there is an one-to-one correspondence between a conic Kropina metric $F$ and a pair $(h,W)$ with $\|W\|_{h}=1$.

In the following, we just study conic Kropina metrics and we always use Kropina metric to take place of conic Kropina metric. From (\ref{eqa4}) or (\ref{NaKro}), we have
\begin{equation}
a_{ij}=e^{-2\rho}h_{ij},\ \ b_{i}=2e^{-2\rho}W_{i},\ \ b^{2}=4e^{-2\rho}, \label{solutionhw1}
\end{equation}
where $\rho := \ln\frac{2}{b}$.

For a Kropina metric  $F=\frac{\alpha ^2}{\beta}$ with navigation data $(h, W)$, let
\begin{equation*}
{\cal R}_{ij}:=\frac{1}{2}(W_{i|j}+W_{j|i}),\ \ {\cal S}_{ij}:=\frac{1}{2}(W_{i|j}-W_{j|i}),
\end{equation*}
\begin{equation*}
{\cal S}^{i}_{\ j}:=h^{ik}\mathcal{S}_{kj},\ \ {\cal S}_{j}:=W^{i}\mathcal{S}_{ij},\ \ \mathcal{R}_{j}:=W^{i}\mathcal{R}_{ij},\ \ \mathcal{R}:=\mathcal{R}_{j}W^{j},
\end{equation*}
where $`` |"$ denotes the covariant derivative with respect to $h$. It is clear that $\mathcal{S}_{j}W^{j}=0$. From (\ref{solutionhw1}), it is easy to get
\begin{equation}\label{rijsij}
\begin{aligned}
&r_{ij}=2e^{-2\rho}(\mathcal{R}_{ij}-W^{l}\rho_{l}h_{ij}),\\
&s_{ij}=2e^{-2\rho}(\mathcal{S}_{ij}+\rho_{i}W_{j}-\rho_{j}W_{i}),
\end{aligned}
\end{equation}
where $\rho _{i}:=\frac{\partial \rho}{\partial x^{i}}$. Then we have
\begin{equation}\label{r00si0s0}
\begin{aligned}
&r_{00}=2e^{-2\rho}(\mathcal{R}_{00}-W^{j}\rho_{j}h^{2}),\\
&{s^{i}}_{0}=2({\mathcal{S}^{i}}_{0}+\rho^{i}W_{0}-\rho_{0}W^{i}),\\
&s_{0}=4e^{-2\rho}(\mathcal{S}_{0}+W^{j}\rho_{j}W_{0}-\rho_{0}),
\end{aligned}
\end{equation}
 where $\rho_{0}:= \rho _{i} y^{i}$ (\cite{ZS}). By conformal relation between $\alpha$ and $h$ given by (\ref{solutionhw1}), the geodesic coefficients $G^{i}_{\alpha}$ of $\alpha$ are related to the geodesic coefficiets $G^{i}_{h}$ of $h$ by
\be
G^{i}_{\alpha}=G^{i}_{h}+\frac{1}{2}\rho^{i}h^{2}-\rho_{0}y^{i},  \label{geodesic3}
\ee
where $\rho ^{i}:=h^{ij}\rho_{j}$.  By (\ref{geodesic2}),(\ref{r00si0s0}) and (\ref{geodesic3}), we obtain the following

\begin{prop} For a Kropina metric $F$ with naviation data $(h, W)$, the geodesic coefficients $G^{i}$ of $F$ can be expressed in terms of the geodesic coefficients $G^{i}_{h}$ of $h$ and the covariant derivatives of $W$ with respect to $h$ as follow:
\be
G^{i}=G^{i}_{h}-F{\mathcal{S}^{i}}_{0}-\frac{1}{2F}(\mathcal{R}_{00}+2F\mathcal{S}_{0})(y^{i}-FW^{i}). \label{geodesic4}
\ee
\end{prop}

\section{S-curvature and Ricci curvature of weakly weighted Einstein-Kropina metrics}\label{Sec3}

In this section, we mainly consider S-curvature and Ricci curvature of weakly weighted Einstein-Kropina metrics. Firstly, we need the following lemma.

\begin{lem}{\rm (\cite{ZS})} For Kropina metric $F=\frac{\alpha^{2}}{\beta}$, we have
\be
\mathbf{S}_{BH}(x,y)=\frac{n+1}{b^{2}}\left(r_{0}-\frac{1}{F}r_{00}\right).  \label{SBH}
\ee
\end{lem}
\vskip 2mm

Now, take $dV=e^{-(n+1)f}dV_{BH}$. Then we have
\beq
&& \mathbf{S}=\mathbf{S}_{BH}+(n+1)f_{0}, \label{S} \\
&& \mathbf{\dot{S}}=\mathbf{\dot{S}}_{BH}+(n+1){\rm Hess}_{F}f(y), \label{Sdot}
\eeq
where $f_{0}:=f_{i}y^{i}$, $f_{i}:=\frac{\pa f}{\pa x^{i}}$ and ${\rm Hess}_{F}f(y)$ is difined by (\ref{Hessian}).

\begin{prop} Let $(M,F,dV =e^{-(n+1)f}dV_{BH})$ be an n-dimensional Kropina measure space. Then
\beq
\frac{1}{n+1}\mathbf{\dot{S}}& =&\frac{1}{b^{2}}\left(r_{0;0}-\frac{\beta}{\alpha^{2}}r_{00;0}+\frac{\alpha^{2}}{\beta}{s^{i}}_{0}r_{i}-2r_{0i}{s^{i}}_{0}\right)
+\frac{1}{b^{4}}\Big[-\big(\frac{\alpha^{2}}{\beta}s_{0}+r_{00}\big)r  \nonumber\\
&& +\frac{2\beta}{\alpha^{2}}r_{00}(3r_{0}-s_{0})+2r_{0}(s_{0}-r_{0})-4\frac{\beta^{2}}{\alpha^{4}}r_{00}^{2}\Big]+ {\rm Hess}_{F}f(y). \label{Sdot1}
\eeq
\end{prop}
\noindent {\it Proof.} \  By (\ref{geodesic2}) and a direct calculation, we have
\begin{equation}
\begin{aligned} \label{Sdot2}
\mathbf{\dot{S}}_{BH}=&\left(\frac{\partial \mathbf{S}_{BH}}{\partial x^{k}}-\frac{\partial G^{j}}{\partial y^{k}}\frac{\partial \mathbf{S}_{BH}}{\partial y^{j}}\right)y^{k}= y^{k}\frac{\partial \mathbf{S}_{BH}}{\partial x^{k}}-2G^{j}\frac{\partial \mathbf{S}_{BH}}{\partial y^{j}}   \\
=&y^{k}\frac{\partial \mathbf{S}_{BH}}{\partial x^{k}}-2G^{j}_{\alpha}\frac{\partial \mathbf{S}_{BH}}{\partial y^{j}}-2T^{j}\frac{\partial \mathbf{S}_{BH}}{\partial y^{j}}  \\
=&(\mathbf{S}_{BH})_{;0}-2T^{j}(\mathbf{S}_{BH})_{y^{j}},
\end{aligned}
\end{equation}
where $T^{j}$ is given by (\ref{geodesicT}) and $`` ; " $ denotes the covariant derivative with respect to $\alpha$. From (\ref{SBH}), we get
\be
(\mathbf{S}_{BH})_{;0}=(n+1)\left[\left(\frac{1}{b^{2}}\right)_{;0}\Big(r_{0}-\frac{1}{F}r_{00}\Big)+\frac{1}{b^{2}}\left(r_{0;0}-\Big(\frac{r_{00}}{F}\Big)_{;0}\right)\right]. \label{SBH;0}
\ee
It is easy to see that
\[
\left(\frac{1}{b^{2}}\right)_{;0}= -\frac{2 (r_{0}+s_{0})}{b^{4}}, \ \  \ \left(\frac{r_{00}}{F}\right)_{;0}= \frac{r_{00}^{2}+\beta r_{00;0}}{\alpha^{2}},
\]
where we have used $\beta_{;0}=r_{00}$. Substituting these into (\ref{SBH;0}), we have
\be
(\mathbf{S}_{BH})_{;0}=(n+1)\left[\frac{1}{b^{2}}\Big(r_{0;0}-\frac{\beta}{\alpha^{2}}r_{00;0}-\frac{r_{00}^{2}}{\alpha^{2}}\Big)+\frac{2}{b^{4}}(r_{0}+s_{0})\Big(\frac{\beta}{\alpha^{2}}r_{00}-r_{0}\Big)\right]. \label{SBH;01}
\ee

On the other hand, it follows from (\ref{SBH}) that
\be
\frac{1}{n+1}(\mathbf{S}_{BH})_{y^{k}}= \frac{1}{b^{2}}\left[r_{k}-\Big(\frac{b_{k}}{\alpha^{2}}-\frac{2\beta y_{k}}{\alpha^{4}}\Big)r_{00}-\frac{2\beta}{\alpha^{2}}r_{0k}\right], \label{SBHyk}
\ee
where $y_{k}=a_{jk}y^{j}$.

Plugging (\ref{geodesicT}),(\ref{SBH;01}) and (\ref{SBHyk}) into (\ref{Sdot2}) yields
\beq
\frac{1}{n+1}\mathbf{\dot{S}}_{BH}& =&\left[\frac{1}{b^{2}}\Big(r_{0;0}-\frac{\beta}{\alpha^{2}}r_{00;0}-\frac{r_{00}^{2}}{\alpha^{2}}\Big)+\frac{2}{b^{4}}(r_{0}+s_{0})\Big(\frac{\beta}{\alpha^{2}}r_{00}-r_{0}\Big)\right] \nonumber\\
&& -\frac{2}{(n+1)}T^{j}(\mathbf{S}_{BH})_{y^{j}}, \label{SBHdot}
\eeq
where
\beq
\frac{2}{n+1}T^{j}(\mathbf{S}_{BH})_{y^{j}}&=&-\left[\frac{\alpha^{2}}{\beta}{s^{j}} _{0}-\frac{1}{b^{2}}\Big(\frac{\alpha^{2}}{\beta}s_{0}+r_{00}\Big)b^{j}+\frac{2}{b^{2}}\Big(s_{0}+\frac{\beta}{\alpha^{2}}r_{00}\Big)y^{j}\right] \nonumber\\
&& \times \frac{1}{b^{2}}\left[r_{j}-\Big(\frac{b_{j}}{\alpha^{2}}-\frac{2\beta y_{j}}{\alpha^{4}}\Big)r_{00}-\frac{2\beta}{\alpha^{2}}r_{0j}\right] \nonumber\\
&=&-(\pi_{1}+\pi_{2}+\pi_{3}) \label{TjSBHyj}
\eeq
and
\beqn
\pi_{1}&:=&\frac{\alpha^{2}}{\beta}{s^{j}}_{0}\times \frac{1}{b^{2}}\left[r_{j}-\Big(\frac{b_{j}}{\alpha^{2}}-\frac{2\beta y_{j}}{\alpha^{4}}\Big)r_{00}-\frac{2\beta}{\alpha^{2}}r_{0j}\right]\nonumber\\
&=& \frac{1}{b^{2}}\left(\frac{\alpha^{2}}{\beta}{s^{j}}_{0}r_{j}-\frac{1}{\beta}s_{0}r_{00}-2r_{0j}{s^{j}}_{0}\right),\\
\pi_{2}&:=&-\frac{1}{b^{2}}\left(\frac{\alpha^{2}}{\beta}s_{0}+r_{00}\right)b^{j}\times \frac{1}{b^{2}}\left[r_{j}-\Big(\frac{b_{j}}{\alpha^{2}}-\frac{2\beta y_{j}}{\alpha^{4}}\Big)r_{00}-\frac{2\beta}{\alpha^{2}}r_{0j}\right]\nonumber\\
&=& \frac{1}{b^{4}}\left[2r_{0}s_{0}-\Big(\frac{\alpha^{2}}{\beta}s_{0}+r_{00}\Big)r-\frac{2\beta}{\alpha^{2}}(s_{0}-r_{0})r_{00}-\frac{2\beta^{2}}{\alpha^{4}}r^{2}_{00}\right]+\frac{1}{b^{2}}\left(\frac{1}{\beta}r_{00}s_{0}+\frac{1}{\alpha^{2}}r^{2}_{00}\right),\\
\pi_{3}&:=&\frac{2}{b^{4}}\left(s_{0}+\frac{\beta}{\alpha^{2}}r_{00}\right)y^{j}\times \left[r_{j}-\Big(\frac{b_{j}}{\alpha^{2}}-\frac{2\beta y_{j}}{\alpha^{4}}\Big)r_{00}-\frac{2\beta}{\alpha^{2}}r_{0j}\right]\nonumber\\
&=& \frac{2}{b^{4}}\left[r_{0}s_{0}+\frac{\beta}{\alpha^{2}}r_{00}(r_{0}-s_{0})-\frac{\beta^{2}}{\alpha^{4}}r^{2}_{00}\right].
\eeqn
Substituting (\ref{TjSBHyj}) into (\ref{SBHdot}) yields
\beq
\frac{1}{n+1}\mathbf{\dot{S}}_{BH}&=&\frac{1}{b^{2}}\left(r_{0;0}-\frac{\beta}{\alpha^{2}}r_{00;0}+\frac{\alpha^{2}}{\beta}{s^{i}}_{0}r_{i}-2r_{0i}{s^{i}}_{0}\right)
+\frac{1}{b^{4}}\Big[-\Big(\frac{\alpha^{2}}{\beta}s_{0}+r_{00}\Big)r \nonumber\\
&& +\frac{2\beta}{\alpha^{2}}r_{00}(3r_{0}-s_{0})+2r_{0}(s_{0}-r_{0})-4\frac{\beta^{2}}{\alpha^{4}}r_{00}^{2}\Big].\label{SBHdot1}
\eeq
By (\ref{Sdot}) and (\ref{SBHdot1}), we get (\ref{Sdot1}).\qed
\vskip 2mm

In order to prove our main theorems, we need the following lemmas.

\begin{lem}{\rm (\cite{XQ})}\label{XiaL} For a Kropina metric $F$ on an n-dimensional manifold $M$, the following are equivalent.
\ben
\item[ {\rm (a)}] $F$ is of isotropic S-curvature with respect to the Busemann-Hausdorff volume form, i.e., $\mathbf{S}_{BH}=(n+1)cF$;
\item[ {\rm (b)}] $r_{00}=\eta \alpha^{2}$;
\item[{\rm (c)}] $\mathbf{S}_{BH}=0$;
\item[{\rm (d)}] $\beta$ is a conformal form with respect to $\alpha$,
\een
where  $\eta = \eta (x)$ is a scalar function on $M$.
\end{lem}

\begin{lem}{\rm (\cite{ZS})}\label{ZSXQ}  For a Kropina metric $F$ on an n-dimensional manifold $M$, $r_{00}=\eta \alpha^{2}$ is equivalent to $\mathcal{R}_{ij}=0$. In this case, $W^{k}\rho_{k}=-\frac{1}{2}\eta$.
\end{lem}

For a Kropina metric $F=\frac{\alpha ^2}{\beta}$ with navigation data $(h,W)$, note that $\|W\|_{h}=1$ and $\frac{\pa}{\pa x^{j}}\|W\|_{h}^{2}= 2({\cal R}_{j}+{\cal S}_{j})=0$. We have the following useful lemma from Lemma \ref{XiaL} and Lemma \ref{ZSXQ}.

\begin{lem}\label{Sj0}  Let  $F=\frac{\alpha ^2}{\beta}$ be a Kropina metric with navigation data $(h,W)$. If $F$ is of isotropic S-curvature with respect to the Busemann-Hausdorff volume form, $\mathbf{S}_{BH}=(n+1)cF$, then ${\cal S}_{j}=0$.

\end{lem}

Now we are in the position to prove our first main result.

\begin{thm} \label{mainth1}  Let $F$ be a  weakly weighted Einstein-Kropina metric on an $n$-dimensional manifold $M$ with volume form $dV=e^{-(n+1)f}dV_{BH}$.  Assume that $\nu \neq 0$. Then $F$ is of isotropic $\mathbf{S}$-curvature with respect to the Busemann-Hausdorff volume form.
\end{thm}
\noindent{\it Proof}. By the assumption, we have
\begin{equation}\label{Ricac}
\begin{aligned}
0=&b^{4}\beta^{2}\alpha^{4}[Ric_{a,c}-(n-1)(3\theta F+\sigma F^{2})]\\
=&b^{4}\beta^{2}\alpha^{4}(Ric+a\mathbf{\dot{S}}-c\mathbf{S}^{2})-(n-1)(3\theta b^{4}\beta\alpha^{6}+\sigma b^{4}\alpha^{8})
\end{aligned}
\end{equation}
for some constants $a$ and $c$. Substituting (\ref{Ricci2}), (\ref{S}) and (\ref{Sdot1}) into (\ref{Ricac}), we have
\begin{equation}\label{Ricac1}
\begin{aligned}
0=&\nu\beta^{4}r_{00}^{2}+\beta^{3}\alpha^{2}\big[\kappa(b^{2}r_{00;0}+2r_{00}r_{0}+2r_{00}s_{0})-2\nu r_{00}r_{0}+2c(n+1)^{2}b^{2}r_{00}f_{0}\big]\\
&+\beta^{2}\alpha^{4}\Big[b^{4}Ric^{\alpha}+b^{2}b^{k}r_{00;k}+(n-2)b^{2}s_{0;0}+b^{2}r_{00}{r^{k}}_{k}-(n-2)s_{0}^{2}\\
&+(-\kappa+n-2)b^{2}r_{0;0}+(\kappa-n)r_{00}r+(-2\kappa+4-2n)r_{0}s_{0}+2(\kappa+1)b^{2}r_{0k}s^{k}_{\ 0}\\
&+(\nu-2\kappa+2-n)r_{0}^{2}+a(n+1)b^{4}{\rm Hess}_{F}f(y)-c(n+1)^{2}(2b^{2}r_{0}f_{0}+b^{4}f_{0}^{2})\Big]\\
&+\beta\alpha^{6}\Big[(\kappa-n)s_{0}r+b^{2}b^{k}s_{0;k}+b^{2}s_{0}r^{k}_{\ k}-b^{4}s^{k}_{\ 0;k}+(-\kappa+n-2)b^{2}r_{k}s^{k}_{\ 0}\\
&-b^{2}r_{0k}s^{k}+(n-1)b^{2}s_{k}s^{k}_{\ 0}-3(n-1)\theta b^{4}\Big]-\alpha^{8}\Big[\frac{1}{2}b^{2}s^{k}s_{k}+\frac{1}{4}b^{4}s^{j}_{\ k}s^{k}_{\ j}+(n-1)\sigma b^{4}\Big],
\end{aligned}
\end{equation}
where we have used that $\kappa =(n-1)-a(n+1)$ and $\nu=3(n-1)-4a(n+1)-c(n+1)^{2}$. (\ref{Ricac1}) can be reorganized as follows:
\be
0=\nu \beta^{4}r_{00}^{2}+\alpha^{2}\beta^{3}P_{1}+\alpha^{4}\beta^{2}P_{2}+\alpha^{6}\beta P_{3}+\alpha^{8}P_{4}. \label{formula3}
\ee

The above equation shows that $\nu\beta^{4} r_{00}^{2}$ can be divided by  $\alpha^{2}$. Since $\beta^{4}$ can not be divided by $\alpha^{2}$ and $\alpha^{2}$ is an irreducible polynomial in $y$, we conclude that there exists a scalar function $\eta (x)$ such that
\[
r_{00}=\eta (x)\alpha^{2}.
\]
Thus, by Lemma \ref{XiaL}, the S-curvature with respect to the Busemann-Hausdorff volume form is isotropic,  $\mathbf{S}_{BH}=0$. \qed

\vskip 2mm

In the following, we will determine the Ricci curvature of weakly weighted Einstein-Kropina metrics.

\begin{lem} Let $h=\sqrt{h_{ij}y^{i}y^{j}}$ be a Riemannian metric on an n-dimensional manifold $M$. Let $W=W^{i}\frac{\partial}{\partial x^{i}}$ be a vector field on $M$ satisfying $\mathcal{R}_{ij}=0$. Then
\be
W_{k|i|j}= -W_{m}\overline{R}^{\ m}_{j \ ki}, \label{formula4}
\ee
\end{lem}
where $\overline{R}^{\ m}_{j \ ki}$ denote the coefficients of the Riemann curvature tensor of $h$.

\noindent{\it Proof}. By the assumption, we have
\be
W_{i|j}+W_{j|i} =0.\label{Rij}
\ee
First, differentiating (\ref{Rij}) and exchanging the indices, we obtain
\beq
&& W_{i|j|k}+W_{j|i|k}=0, \label{diffRij1} \\
&& W_{j|k|i}+W_{k|j|i}=0, \label{diffRij2} \\
&& W_{k|i|j}+W_{i|k|j}=0. \label{diffRij3}
\eeq
Adding (\ref{diffRij2}) and (\ref{diffRij3}) together, then the sum being subtracted by (\ref{diffRij1}) , we obtain
\be
(W_{i|k|j}-W_{i|j|k})+(W_{j|k|i}-W_{j|i|k})+(W_{k|j|i}-W_{k|i|j})+2W_{k|i|j}=0. \label{formula5}
\ee
Using the Ricci identity, $W_{k|i|j}-W_{k|j|i}=W_{m}\overline{R}^{~m}_{k~ij}$, we obtain
\be
W_{m}\overline{R}^{~m}_{i~kj}+W_{m}\overline{R}^{~m}_{j~ki}+W_{m}\overline{R}^{~m}_{k~ji}+2W_{k|i|j}=0. \label{formula6}
\ee
By applying the Bianchi identity, $\overline{R}^{~m}_{i~kj}+\overline{R}^{~m}_{k~ji}+\overline{R}^{~m}_{j~ik}=0$,
we obtain (\ref{formula4}).  \qed

\vskip 2mm

By (\ref{geodesic4}), we have
\[
G^{i}=G^{i}_{h}+Q^{i},
\]
where
\[
Q^{i}:=-F{\mathcal{S}^{i}}_{0}-\frac{1}{2F}(\mathcal{R}_{00}+2F\mathcal{S}_{0})(y^{i}-FW^{i}).
\]
Then, by (\ref{Riemann1}), we have
\be
R^{i}_{\ k}=\overline{R}^{i}_{\ k}+2{Q^{i}}_{|k}-[{Q^{i}}_{|m}]_{y^{k}}y^{m}+2Q^{m}[Q^{i}]_{y^{m}y^{k}}-[Q^{i}]_{y^{m}}[Q^{m}]_{y^{k}}.\label{Riemann2}
\ee
Here $ `` |"$ denotes the covariant differentiation with respect to $h$.

From now on, we assume that $F$ is of isotropic S-curvature with respect to the Busemann-Hausdorff volume form. By Lemma \ref{XiaL}, Lemma \ref{ZSXQ} and Lemma \ref{Sj0}, $\mathcal{R}_{ij}=0$ and ${\cal S}_{j}=0$. Then the geodesic coefficients $G^{i}$ are reduced to the following expression:
\be
G^{i}=G^{i}_{h}+Q^{i},\label{geodesic5}
\ee
where
\be
Q^{i}=-F\mathcal{S}^{i}_{~0}.\label{Qi}
\ee
By (\ref{formula4}), we have
\be
W_{i|j|k}=-W^{p}{\overline{R}}_{kpij}. \label{formula7}
\ee
Then we have
\beq
&& {\mathcal{S}^{i}}_{0|k}=\overline{R}^{~i}_{p~kq}y^{p}W^{q}, \label{S1}\\
&& {\mathcal{S}^{i}}_{0|0}=-\overline{R}^{~i}_{p~mq}y^{p}y^{q}W^{m}, \label{S2}\\
&& {\mathcal{S}}_{0|k}=\mathcal{S}_{mk}\mathcal{S}^{m}_{\ 0}+W_{m}\overline{R}^{\ m}_{p~kq}y^{p}W^{q}. \label{S3}
\eeq
Observe that
\[
W_{0|k}=\mathcal{S}_{0k}, \ \ {W^{i}}_{|k}={\mathcal{S}^{i}}_{k}, \ \ F_{|k}=\frac{F}{W_{0}}\mathcal{S}_{k0}.
\]
Now, for simplicity, let
\be
 \xi^{i}:=y^{i}-FW^{i}. \label{Aiji}
\ee
By a direct calculation, we have
\be
{Q^{i}}_{|k}=-F\overline{R}^{~i}_{p~kq}y^{p}W^{q}-F_{|k}{\cal S}^{i}_{\ 0}. \label{Qi|k}
\ee

Futher, we need the following formula,
\[
F_{y^{k}}=\frac{y_{k}-FW_{k}}{W_{0}}=\frac{\xi_{k}}{W_{0}},
\]
where $y_{k}=h_{kj}y^{j}$ and $\xi_{k}:=h_{ik}\xi^{i}$.

Thus, by a series computations, we obtain that
\beq
{[{Q^{i}}_{|m}]}_{y^{k}}y^{m}&=&\frac{\xi_{k}}{W_{0}}\overline{R}^{~i}_{p~mq}y^{p}y^{q}W^{m}-F\overline{R}^{~i}_{k~mq}W^{q}y^{m}+F_{|k}{\cal S}^{i}_{\ 0}. \\
Q^{m}[Q^{i}]_{y^{m}y^{k}}&=&\frac{F}{W_{0}}\left({\cal S}^{i}_{\ 0}\mathcal{S}_{k0}+ \xi_{k}{\mathcal{S}^{m}}_{0}{\cal S}^{i}_{\ m}\right).\\
{[Q^{i}]_{y^{m}}}{[Q^{m}]_{y^{k}}}&=& F^{2}{\cal S}^{m}_{\ k}{\cal S}^{i}_{\ m}+\frac{F}{W_{0}}\left(- {\cal S}_{k0}{\cal S}^{i}_{\ 0} +\xi_{k}{\mathcal{S}^{i}}_{m}{\mathcal{S}^{m}}_{0}\right). \label{QimQmk}
\eeq
Plugging (\ref{Qi|k})-(\ref{QimQmk}) into (\ref{Riemann2}) yields the following result.

\begin{prop} Let $F$ be a Kropina metric expressed by (\ref{hwKropina}). Suppose that it is of isotropic S-curvature with respect to the Busemann-Hausdorff volume form. Then the Riemann curvature of $F$ can be expressed in terms of the Riemann curvature of $h$ and the covariant derivatives of $W$ with respect to $h$ as follows:
\begin{align}
{R^{i}}_{k}&={\overline{R}^{i}}_{k}-2F\overline{R}^{\ i}_{p \ kq}y^{p}W^{q}-\frac{\xi_{k}}{W_{0}}\overline{R}^{\ i}_{p \ mq}y^{p}y^{q}W^{m}+F\overline{R}^{\ i}_{k \ mq}y^{m}W^{q}\nonumber\\
& -F^{2}\mathcal{S}^{m}_{\ k}\mathcal{S}^{i}_{\ m} +\frac{\xi_{k}}{W_{0}}  F\mathcal{S}^{m}_{\ 0}{\cal S}^{i}_{\ m}.\label{RiemanNG}
\end{align}
\end{prop}

From (\ref{RiemanNG}), we have the following

\begin{prop}\label{Riccibyh} Let $F$ be a Kropina metric expressed by (\ref{hwKropina}). Suppose that it is of isotropic S-curvature with respect to the Busemann-Hausdorff volume form.  Then the Ricci curvature of $F$ can be expressed in terms of the Ricci curvature of $h$ and the covariant derivatives of $W$ with respect to $h$ as follow:
\beq
Ric &=&\overline{Ric}- 2F\overline{R}^{\ i}_{p \ iq}y^{p}W^{q}-F^{2}{\cal S}^{m}_{\ i}{\cal S}^{i}_{\ m},\label{RicNG}
\eeq
where $\overline{Ric}$ denotes the Ricci curvature of $h$.
\end{prop}
\noindent{\it Proof}. By applying the Bianchi identity, $\overline{R}^{\ i}_{i \ mq}+\overline{R}^{\ i}_{m \ qi}+\overline{R}^{\ i}_{q \ im}=0$, we have
\beq
F\overline{R}^{\ i}_{i \ mq}y^{m}W^{q}&=&F\big(-\overline{R}^{\ i}_{m \ qi}-\overline{R}^{\ i}_{q \ im}\big)y^{m}W^{q} \nonumber\\
&=&F\big(\overline{R}^{\ i}_{m \ iq}-\overline{R}^{\ i}_{q \ im}\big)y^{m}W^{q}=0.
\eeq
By (\ref{S3}) and ${\cal S}_{j}=0$, we have ${\cal S}^{m}_{\ 0}{\cal S}^{i}_{\ m}=W^{m}\overline{R}^{\ i}_{p \ mq}W^{p}y^{q}$. Then
\beqn
&& \frac{\xi_{i}}{W_{0}} F\mathcal{S}^{m}_{\ 0}{\cal S}^{i}_{\ m}-\frac{\xi_{i}}{W_{0}}\overline{R}^{\ i}_{p \ mq}y^{p}y^{q}W^{m}=\frac{\xi_{i}}{W_{0}}\overline{R}^{\ i}_{p \ mq}W^{m}y^{q}(FW^{p}-y^{p})\\
&&=- \frac{\xi_{i}}{W_{0}}\overline{R}^{\ i}_{p \ mq} \xi ^{p}W^{m}y^{q}=0.
\eeqn
Then (\ref{RicNG}) follows from (\ref{RiemanNG}). \qed

\vskip 2mm

\section{The weakly weighted Einstein-Kropina metrics with $\nu\neq 0$}\label{Sec4}

In this section, we will firstly characterize weakly weighted Einstein-Kropina metrics via navigation data $(h, W)$ in the case that $\nu\neq 0$.

\begin{thm} \label{mainth3}
Let $F$ be a Kropina metric on an n-dimensional manifold $M$ defined by a navigation data $(h,W)$.  Then $F$ is a weakly weighted Einstein metric with weight constants $a$ and $c$ satisfying
\be
{\rm Ric}_{a,c}=(n-1)\left(\frac{3\theta}{F}+\sigma\right)F^{2} \label{Ricac2}
\ee
with respect to a volume form $dV=e^{-(n+1)f}dV_{BH}$ if and only if  there exists a scalar function $\mu$ on $M$, such that the Ricci curvature tensor of $h$ satisfies
\be
 Ric^{h}+a(n+1){\rm Hess}_{h}f - c(n+1)^{2}(df \otimes df)= (n-1)\mu {h}^{2} \label{formula13}
\ee
and $W$ satisfies ${\cal R}_{ij}=0$. In this case,
\be
\sigma = \mu-\frac{1}{n-1}\big\{Ric^{h}(W)+\mathcal{S}^{p}_{\ q}\mathcal{S}^{q}_{\ p}+a(n+1){\rm Hess}_{h}f(W)-c(n+1)^{2}(f_{p}W^{p})^{2}\big\} \label{sigma2}
\ee
and
\be
\theta_{i}= \frac{1}{3(n-1)}\left[2a(n+1)(f_{ip}W^{p}+f_{p}\mathcal{S}^{p}_{\ i})-2c(n+1)^{2}f_{i}f_{p}W^{p}\right].\label{thetap}
\ee
\end{thm}

\noindent{\it Proof}. \ Firstly, suppose that $F$ is a weakly weighted Einstein metric satisfying (\ref{Ricac2}).  By Theorem \ref{mainth1},  $\mathbf{S}_{BH}=0$, that is, $\mathcal{R}_{ij}=0$ and ${\cal S}_{j}=0$. Then by (\ref{geodesic4}), we have
\be
G^{i}_{h}-G^{i}=F\mathcal{S}^{i}_{\ 0}. \label{geodesic6}
\ee
Recall that $F(x,y)=h(x,\xi):=\tilde{h}$ and $\xi:=y-FW$. We have
\begin{equation}
G^{i}_{h}-G^{i}=(\xi^{j}+\tilde{h}W^{j})\tilde{h}\mathcal{S}^{i}_{~j}.\label{geodesic7}
\end{equation}
In this case, ${\dot{\bf S}}= \dot{\bf S}_{BH}+(n+1){\rm Hess}_{F}f(y) =(n+1){\rm Hess}_{F}f(y)$. Further, by (\ref{Hessian}), we have
\beq
{\rm Hess}_{F}f(y) &= &f_{ij}(\xi^{i}+\tilde{h}W^{i})(\xi^{j}+\tilde{h}W^{j})+2f_{i}(G^{i}_{h}-G^{i})\nonumber\\
&=& f_{ij}(\xi^{i}+\tilde{h}W^{i})(\xi^{j}+\tilde{h}W^{j})+2f_{i}(\xi^{j}+\tilde{h}W^{j})\tilde{h}\mathcal{S}^{i}_{~j}. \label{HessF}
\eeq
Here, $f_{i}:=\frac{\partial f}{\partial x^{i}}$ and $ f_{ij}:=\frac{\partial^{2}f}{\partial x^{i}\partial x^{j}}-\overline{\Gamma}^{m}_{ij}\frac{\partial f}{\partial x^{m}}$ denote the second order covariant derivatives of $f$ with respect to $h$.

By Proposition \ref{Riccibyh}, (\ref{Ricac2}) is equivalent to
\begin{align}\label{Ricac3}
&Ric^{h}(\xi)-\tilde{h}^{2}Ric^{h}(W)-\tilde{h}^{2} \mathcal{S}^{m}_{~i}\mathcal{S}^{i}_{~m}
+a(n+1)\big[f_{pq}(\xi^{p}\xi^{q}+2\tilde{h}\xi^{p}W^{q} \nonumber\\
&+\tilde{h}^{2}W^{p}W^{q})+2\tilde{h}f_{i}\mathcal{S}^{i}_{~p}\xi^{p}\big] -c(n+1)^{2}\big[\tilde{f}_{0}^{2}+2\tilde{h}\tilde{f}_{0}f_{p}W^{p}+\tilde{h}^{2}(f_{p}W^{p})^{2}\big] \nonumber\\
&=(n-1)[3\tilde{\theta} \tilde{h}+(3\theta_{p}W^{p}+\sigma)\tilde{h}^{2}],
\end{align}
where $\widetilde{W}_{0}:= W_{p}\xi^{p}, \ \widetilde{\cal S}_{0}:={\cal S}_{p}\xi ^{p}, \ \tilde{f}_{0}:= f_{p}\xi ^{p}$ and $\tilde{\theta}:=\theta_{i}\xi^{i}$.

Note that $\tilde{h}$ is irrational in $\xi$. Separating rational and irrational terms in the above equation we have
\beq
&& Ric^{h}(\xi)-\tilde{h}^{2}Ric^{h}(W) -\tilde{h}^{2} \mathcal{S}^{m}_{\ i}\mathcal{S}^{i}_{\ m} +a(n+1)\big[{\rm Hess}_{h}f(\xi)+\tilde{h}^{2}{\rm Hess}_{h}f(W) \big] \nonumber\\
&& -c(n+1)^{2}\big[\tilde{f}_{0}^{2}+\tilde{h}^{2}(f_{p}W^{p})^{2}\big] =(n-1)(3\theta_{p}W^{p}+\sigma)\tilde{h}^{2} \label{rational}
\eeq
and
\be
2a(n+1)\xi^{p}(f_{pq}W^{q}+f_{q}\mathcal{S}^{q}_{\ p})-2c(n+1)^{2}\tilde{f}_{0}f_{p}W^{p}=3(n-1)\tilde{\theta}.\label{irrational}
\ee

\vskip 2mm

From (\ref{rational}) we obtain
\beq
Ric^{h}(\xi)+a(n+1){\rm Hess}_{h}f(\xi)-c(n+1)^{2}\tilde{f}_{0}^{2}=(n-1)\mu\tilde{h}^{2}, \label{formula14}
\eeq
where
\be
\mu = 3\theta_{p}W^{p}+\sigma+\frac{1}{n-1}\big\{Ric^{h}(W) +\mathcal{S}^{m}_{\ i}\mathcal{S}^{i}_{\ m} -a(n+1){\rm Hess}_{h}f(W)+c(n+1)^{2}(f_{p}W^{p})^{2}\big\}.
\ee
From (\ref{formula14}), we get (\ref{formula13}).

Conversely, suppose that (\ref{formula14}) holds and $\mathcal{R}_{ij}=0$.  We choose
\be
\theta_{i}:= \frac{1}{3(n-1)}\left[2a(n+1)(f_{ip}W^{p}+f_{p}\mathcal{S}^{p}_{\ i})-2c(n+1)^{2}f_{i}f_{p}W^{p}\right]. \label{thetap2}
\ee
And then, equation (\ref{irrational}) holds. Furthermore, we choose
\be
\sigma :=\mu-\frac{1}{n-1}\big\{Ric^{h}(W)+\mathcal{S}^{p}_{\ q}\mathcal{S}^{q}_{\ p}+a(n+1){\rm Hess}_{h}f(W)-c(n+1)^{2}(f_{p}W^{p})^{2}\big\}. \label{sigma3}
\ee
It is easy to check that equation (\ref{rational}) also holds and it follows that
\[
Ric_{a,c}=(n-1)\left(\frac{3\theta}{F}+\sigma\right)F^{2}
\]
with respect to $dV=e^{-(n+1)f}dV_{BH}$. \qed

\vskip 2mm

In the following,  we will characterize weakly weighted Einstein-Kropina metrics via $\alpha$ and $\beta$ in the case that $\nu\neq 0$. Assume that $F=\frac{\alpha^{2}}{\beta}$ is a weakly weighted Einstein-Kropina metric satisfying (\ref{weweEF}) and $\nu\neq0$.  By Theorem \ref{mainth1} and Lemma \ref{XiaL}, we have
\be
r_{00}=\eta(x)\alpha^{2}. \label{isotropic2}
\ee
Then, it is easy to get
\beq
&& r_{0i}=\eta y_{i},\ \ r_{i}=\eta b_{i},\ \ r=\eta b^{2},\ \ r^{i}_{~j}=\eta\delta^{i}_{~j},\nonumber\\
&& r_{0k}s^{k}_{\ 0}=0,\ \ r_{0k}s^{k}=\eta s_{0},\ \ r_{0}=\eta\beta,\ \ s^{k}_{\ 0}r_{k}=\eta s_{0}, \nonumber\\
&& r_{00;k}=\eta_{k}\alpha^{2},\ \ r_{00;0}=\eta_{0}\alpha^{2},\ \ r_{0;0}=\eta_{0}\beta+\eta^{2}\alpha^{2}.\label{funfor}
\eeq
Here, $y_{i}:=a_{ij}y^{j}$. Substituting (\ref{funfor}) into (\ref{Ricac1}) and dividing both sides of (\ref{Ricac1}) by  $\alpha^{4}$, we can obtain
\begin{align}
0=&Ric^{\alpha}b^{4}\beta^{2}+(n-2)\beta^{2}\big[b^{2}(s_{0;0}+\eta_{0}\beta)-2\eta\beta s_{0}-s_{0}^{2}-\eta^{2}\beta^{2}\big] \nonumber\\
&-\left[3\kappa-\nu-a(n+1)\right]b^{4}\beta^{2}f_{0}^{2}+(-\kappa+n-1)b^{4}\beta^{2}{\rm Hess}_{F}f(y) \nonumber\\
&+b^{2}\beta\alpha^{2}\big[\beta b^{k}\eta_{k}+(n-2)\eta^{2}\beta +(n-3)\eta s_{0}+b^{k}s_{0;k}-b^{2}s^{k}_{~0;k}
 \nonumber\\
&+(n-1)(s_{k}s^{k}_{~0}-3b^{2}\theta)\big]-b^{2}\alpha^{4}\Big(\frac{1}{2}s^{k}s_{k}+\frac{b^{2}}{4}s^{j}_{~k}s^{k}_{~j}+(n-1)b^{2}\sigma\Big). \label{Ricac12}
\end{align}
From (\ref{Ricac12}), it is easy to see that there exists some scalar function $ \lambda = \lambda(x)$ on $M$ such that
\begin{align}
&Ric^{\alpha}b^{4}+(n-2)[b^{2}(s_{0;0}+\eta_{0}\beta)-2\eta\beta s_{0}-s_{0}^{2}-\eta^{2}\beta^{2}]\nonumber\\
&-\left(3\kappa-\nu-a(n+1)\right)b^{4}f_{0}^{2}+(-\kappa+n-1)b^{4}{\rm Hess}_{F}f(y)= \lambda \alpha^{2}.\label{Ricablam}
\end{align}
Then (\ref{Ricac12}) can be simplified as
\beq
0&=&\beta\Big\{\lambda\beta+b^{2}[\beta b^{k}\eta_{k}+(n-2)\eta^{2}\beta+(n-3)\eta s_{0}+b^{k}s_{0;k}-b^{2}s^{k}_{ \ 0;k}\nonumber\\
&& +(n-1)(s_{k}s^{k}_{~0}-3b^{2}\theta)]\Big\}-b^{2}\alpha^{2}\Big(\frac{1}{2}s^{k}s_{k}+\frac{b^{2}}{4}s^{j}_{~k}s^{k}_{~j}+(n-1)b^{2}\sigma\Big). \label{Ricac13}
\eeq
Since $\alpha^{2}$ is irreducible polynomial in $y$,  (\ref{Ricac13}) implies the following equations
\be
\lambda \beta+b^{2}[(n-3)\eta s_{0}+(n-2)\eta^{2}\beta+b^{k}\eta_{k}\beta+b^{k}s_{0;k}-b^{2}s^{k}_{~0;k}+(n-1)(s_{k}s^{k}_{~0}-3b^{2}\theta)]=0, \label{formula8}
\ee
\be
\frac{1}{2}s^{k}s_{k}+\frac{b^{2}}{4}s^{j}_{~k}s^{k}_{~j}+(n-1)b^{2}\sigma =0. \label{formula82}
\ee
Differentiating both sides of (\ref{formula8}) with respect to $y^{i}$ yields
\beq
&& \lambda b_{i}+[(n-2)\eta^{2}+b^{k}\eta_{k}]b^{2}b_{i}-3(n-1)b^{4}\theta_{i}+(n-1)b^{2}s_{k}s^{k}_{\ i}+(n-3)\eta b^{2}s_{i} \nonumber\\
&& +b^{2}b^{k}s_{i;k}-b^{4}s^{k}_{\ i;k}=0. \label{formula9}
\eeq
Contracting (\ref{formula9}) with $b^{i}$ gives
\be
\lambda b^{2}+[(n-2)\eta^{2}+b^{k}\eta_{k}]b^{4}-3(n-1)b^{4}\theta_{i}b^{i}-(n-2)b^{2}s_{k}s^{k}+b^{4}(s^{k}_{~;k}+s^{k}_{\ i}s^{i}_{\ k})=0,\label{formula10}
\ee
where we used
\beqn
&& s_{i;k}b^{i}=-\eta s_{k}-s_{i}s^{i}_{\ k},\\
&& s^{k}_{\ i;k}b^{i}=-s^{k}_{~;k}-s^{k}_{\ i}s^{i}_{\ k}.
\eeqn
Hence, we obtain
\be
\lambda =-\left[(n-2)\eta^{2}+b^{k}\eta_{k}\right]b^{2}+3(n-1)b^{2}\theta_{i}b^{i}+(n-2)s_{k}s^{k}-b^{2}(s^{k}_{~;k}+s^{k}_{~i}s^{i}_{~k}).\label{lambda}
\ee
Plugging (\ref{lambda}) into  (\ref{formula8}) yields
\begin{equation}\label{formula12}
\begin{aligned}
0&=\beta[(n-2)s^{k}s_{k}+3(n-1)b^{2}\theta_{i}b^{i}-b^{2}(s^{k}_{~;k}+s^{k}_{~i}s^{i}_{~k})]+b^{2}[(n-3)\eta s_{0}\\
&+b^{k}s_{0;k}-b^{2}s^{k}_{\ 0;k}+(n-1)s_{k}s^{k}_{\ 0}-3(n-1)b^{2}\theta].
\end{aligned}
\end{equation}
Further, by (\ref{formula82}), we obtain
\be
\sigma=-\frac{1}{(n-1)b^{2}}\left(\frac{1}{2}s^{k}s_{k}+\frac{b^{2}}{4}s^{j}_{~k}s^{k}_{~j}\right).\label{sigma}
\ee

Conversely, if (\ref{isotropic2}),(\ref{Ricablam}) and (\ref{formula12}) hold and $\lambda$, $\sigma$ are given by (\ref{lambda}) and (\ref{sigma}) respectively,  then it is easy to see that (\ref{Ricac1}) holds, that is, $F$ is a weakly weighted Einstein-Kropina metric with weight constants $a$ and $c$. Hence, we have actually proved the following

\begin{thm} \label{mainth2}
Let $a,c$ be two constants satisfying $\nu\neq 0$ and $F=\frac{\alpha^{2}}{\beta}$ be a Kropina metrics on an n-dimensional manifold $M$. Then $F$ is a weakly weighted Einstein metric satisfying
\[
Ric_{a,c}=(n-1)\left(\frac{3\theta}{F}+\sigma\right)F^{2}
\]
with respet to some volume form $dV=e^{-(n+1)f}dV_{BH}$ if and only if equations (\ref{isotropic2}), (\ref{Ricablam}) and (\ref{formula12}) are satisfied for some scalar functions $\lambda$, $\sigma$ given by (\ref{lambda}) and (\ref{sigma}) respectively.
\end{thm}

\section{The weakly weighted Einstein-Kropina metrics with $\nu=0$ and $\kappa\neq 0$}\label{Sec5}

By the definition, when $\nu=0$ and $\kappa\neq 0$,  the generalized weighted Ricci curvature with weight constants $a$ and $c$ is given by
\[
{\rm Ric}_{a,c}={\rm PRic}-\kappa\left[\frac{\mathbf{\dot{S}}}{n+1}+\frac{4\mathbf{S}^{2}}{(n+1)^{2}}\right].
\]
In the following, we are going to derive an equivalent condition for a Kropina metric $F$ to satisfy
\be
Ric_{a,c}=(n-1)\left(\frac{3\theta}{F}+\sigma\right)F^{2}. \label{Ricac4}
\ee

Assume that $F$ is a weakly weighted Einstein-Kropina metric satisfying (\ref{Ricac4}). In this case, the equation (\ref{formula3}) becomes
\be
\beta^{3}P_{1}+\beta^{2}\alpha^{2}P_{2}+\beta\alpha^{4}P_{3}+\alpha^{6}P_{4}=0 \label{formula15}
\ee
and the polynomials $P_{i}$'s can be simplified a little:
\beqn
P_{1}&=&\kappa(b^{2}r_{00;0}+2r_{00}r_{0}+2r_{00}s_{0})+2[3\kappa-a(n+1)]b^{2}r_{00}f_{0}, \\
P_{2}&=& b^{4}{\rm Ric}^{\alpha}+b^{2}b^{k}r_{00;k}+(n-2)b^{2}s_{0;0}+b^{2}r_{00}r^{k}_{\ k}-(n-2)s_{0}^{2}+(-\kappa+n-2)b^{2}r_{0;0} \\
&&+(\kappa-n)r_{00}r+(-2\kappa+4-2n)r_{0}s_{0}+2(\kappa+1)b^{2}r_{0k}s^{k}_{\ 0}+(-2\kappa+2-n)r_{0}^{2} \\
&&+(-\kappa+n-1)b^{4}{\rm Hess}_{F}f(y)-(3\kappa-a(n+1))\left(2b^{2}r_{0}f_{0}+b^{4}f_{0}^{2}\right), \\
P_{3}&=&(\kappa-n)s_{0}r+b^{2}b^{k}s_{0;k}+b^{2}s_{0}r^{k}_{\ k}-b^{4}s^{k}_{\ 0;k}+(-\kappa+n-2)b^{2}r_{k}s^{k}_{\ 0} \\
&&-b^{2}r_{0k}s^{k}+(n-1)b^{2}s_{k}s^{k}_{\ 0}-3(n-1)\theta b^{4},\nonumber\\
P_{4}&=& -b^{2}\left[\frac{1}{2}s^{k}s_{k}+\frac{1}{4}b^{2}{s^{j}}_{k}{s^{k}}_{j}+(n-1)\sigma b^{2}\right].
\eeqn
By (\ref{formula15}), we know that there exists a 1-form $ \zeta = \zeta _{i}(x)y^{i}$ such that
\[
P_{1}= \zeta \alpha^{2},
\]
which is equivalent to
\be
\zeta \alpha^{2}= \kappa(b^{2}r_{00;0}+2r_{00}r_{0}+2r_{00}s_{0})+2[3\kappa-a(n+1)]b^{2}r_{00}f_{0}.   \label{zeta}
\ee
Then (\ref{formula15}) can be simplified as
\be
\beta^{2}(\beta \zeta +P_{2})+\beta\alpha^{2}P_{3}+\alpha^{4}P_{4}=0. \label{formula16}
\ee
By (\ref{formula16}), we know there exists some function $u = u (x)$ such that
\be
\beta \zeta+P_{2}= u \alpha^{2}. \label{f2u}
\ee
Then (\ref{formula16}) can be simplified as
\[
\beta(\beta u+P_{3})+\alpha^{2}P_{4}=0,
\]
that is
\beqn
&&\beta\big[\beta u+(\kappa-n)s_{0}r+b^{2}b^{k}s_{0;k}+b^{2}s_{0}{r^{k}}_{k}-b^{4}{s^{k}}_{0;k}+(-\kappa+n-2)b^{2}r_{k}{s^{k}}_{0}
-b^{2}r_{0k}s^{k}\nonumber\\
&&+(n-1)b^{2}s_{k}{s^{k}}_{0}-3(n-1)\theta b^{4}\big]-\alpha^{2}b^{2}\left[\frac{1}{2}s^{k}s_{k}+\frac{1}{4}b^{2}{s^{j}}_{k}{s^{k}}_{j}+(n-1)\sigma b^{2}\right]=0. \label{formula18}
\eeqn
Since $\alpha^{2}$ can't be divided by $\beta$, we see that above equation is equivalent to the following equations
\beq
&&\beta u+(\kappa-n)s_{0}r+b^{2}b^{k}s_{0;k}+b^{2}s_{0}{r^{k}}_{k}-b^{4}{s^{k}}_{0;k}+(-\kappa+n-2)b^{2}r_{k}{s^{k}}_{0}
-b^{2}r_{0k}s^{k}\nonumber\\
&&+(n-1)b^{2}s_{k}{s^{k}}_{0}-3(n-1)\theta b^{4}=0,  \label{formula19} \\
&&\frac{1}{2}s^{k}s_{k}+\frac{b^{2}}{4}s^{j}_{\ k}s^{k}_{\ j}+(n-1)b^{2}\sigma =0.\label{formula192}
\eeq

Differentiating both sides of (\ref{formula19}) with respect to $y^{i}$ yields
\beq
&& u b_{i}+(\kappa-n)s_{i}r+b^{2}s_{i}{r^{k}}_{k}+(-\kappa+n-2)b^{2}r_{k}{s^{k}}_{i}-b^{2}r_{ik}s^{k}+(n-1)b^{2}s_{k}{s^{k}}_{i} \nonumber\\
&&-3(n-1)\theta_{i} b^{4}+b^{2}b^{k}s_{i;k}-b^{4}{s^{k}}_{i;k}=0. \label{formula20}
\eeq
Further, contracting (\ref{formula20}) with $b^{i}$ gives
\be
u b^{2}+(\kappa -n)b^{2}r_{i}s^{i}-(n-2)b^{2}s_{i}s^{i}-3(n-1)\theta_{i}b^{i} b^{4}+b^{4}(s^{k}_{~;k}+s^{k}_{~i}s^{i}_{~k}+s^{k}_{~i}r^{i}_{~k})=0,\label{formula21}
\ee
where we have used
\beqn
&& s_{i;k}b^{i}=-s_{i}(r^{i}_{~k}+s^{i}_{~k}),\\
&& s^{k}_{~i;k}b^{i}=-(s^{k}_{~;k}+s^{k}_{~i}s^{i}_{~k}+s^{k}_{~i}r^{i}_{~k}).
\eeqn
Then we obtain
\be
u =(n-\kappa)r_{i}s^{i}+(n-2)s^{i}s_{i}+3(n-1)b^{2}\theta_{i}b^{i}-b^{2}(s^{k}_{\ ;k}+s^{k}_{\ i}s^{i}_{\ k}+s^{k}_{\ i}r^{i}_{\ k}).\label{f2define}
\ee
Plugging (\ref{f2define}) into (\ref{formula19})  yields
\beq
&&\beta\left[(n-\kappa)r_{i}s^{i}+(n-2)s^{i}s_{i}+3(n-1)b^{2}\theta_{i}b^{i}-b^{2}(s^{k}_{~;k}+s^{k}_{\ i}s^{i}_{\ k}+s^{k}_{\ i}r^{i}_{\ k})\right] \nonumber\\
&&+ (\kappa -n)s_{0}r +b^{2}b^{k}s_{0;k}+b^{2}s_{0}r^{k}_{\ k}-b^{4}s^{k}_{\ 0;k}+(n-\kappa -2)b^{2}r_{k}s^{k}_{\ 0}\nonumber\\
&& -b^{2}r_{0k}s^{k}+(n-1)b^{2}s_{k}{s^{k}}_{0}-3(n-1)\theta b^{4}=0.\label{formula23}
\eeq

Furthermore, from (\ref{formula192}), we obtain
\be
\sigma =-\frac{1}{(n-1)b^{2}}\left(\frac{1}{2}s^{k}s_{k}+\frac{b^{2}}{4}s^{j}_{~k}s^{k}_{~j}\right). \label{formula24}
\ee

Conversely, if   (\ref{zeta}) and (\ref{f2u}), (\ref{formula23}) hold for some 1-form $\zeta$  with (\ref{f2define}) and (\ref{formula24}). Then it is easy to see that $F$ is weakly weighted Einstein-Kropina metric.
Thus we have proved the following

\begin{thm}\label{mainth4}
Let $a,c$  be two constants satisfying $\nu=0$ and $\kappa\neq0$. Let $F=\frac{\alpha^{2}}{\beta}$ be a Kropina metric on an n-dimensional manifold $M$. Then F is  weakly weighted Einstein-Kropina metric satisfying
\[
Ric_{a,c}=(n-1)\left(\frac{3\theta}{F}+\sigma\right)F^{2}
\]
with respect to some volume form $dV=e^{-(n+1)f}dV_{BH}$ if and only if equations  (\ref{zeta}) and (\ref{f2u}), (\ref{formula23})  hold for some 1-form $\zeta$ and $u= u(x)$ and $\sigma$ are determined by (\ref{f2define}) and (\ref{formula24}) respectively.
\end{thm}

\section{The weakly weighted Einstein-Kropina metrics with $\nu=0$ and $\kappa = 0$}\label{Sec6}

In this section we shall consider the weakly weighted Einstein-Kropina metrics with $\nu=0$ and $\kappa=0$. In this case, the generalized weighted Ricci curvatures are just the projective Ricci curvature (\cite{CSM}\cite{SS})
\[
{\rm Ric}_{a,c}={\rm PRic}
\]
and $a=\frac{n-1}{n+1}, \ c= -\frac{n-1}{(n+1)^2}$. Then, the polynomials $P_{i}$'s in (\ref{formula15}) can be simplified as follows.
\beqn
P_{1}&=&-2(n-1)b^{2}r_{00}f_{0}, \\
P_{2}&=&b^{4}Ric^{\alpha}+b^{2}b^{k}r_{00;k}+(n-2)(b^{2}s_{0;0}-s_{0}^{2}+b^{2}r_{0;0}-2r_{0}s_{0}-r_{0}^{2})+b^{2}r_{00}{r^{k}}_{k} \\
&&-nr_{00}r+2b^{2}r_{0k}{s^{k}}_{0}+(n-1)b^{4}{\rm Hess}_{F}f(y)+ (n-1)\left(2b^{2}r_{0}f_{0}+b^{4}f_{0}^{2}\right), \\
P_{3}&=&-ns_{0}r+b^{2}b^{k}s_{0;k}+b^{2}s_{0}{r^{k}}_{k}-b^{4}{s^{k}}_{0;k}+(n-2)b^{2}r_{k}{s^{k}}_{0} \\
&& -b^{2}r_{0k}s^{k}+(n-1)b^{2}s_{k}{s^{k}}_{0}-3(n-1)\theta b^{4}, \\
P_{4}&=&-b^{2}\left[\frac{1}{2}s^{k}s_{k}+\frac{1}{4}b^{2}{s^{j}}_{k}{s^{k}}_{j}+(n-1)\sigma b^{2}\right].
\eeqn
At the same time, (\ref{zeta}) is reduced equivalently to
\be
\zeta \alpha^{2} =-2(n-1)b^{2}r_{00} f_{0}. \label{zeta6}
\ee
From (\ref{zeta6}), we find that there exists a scalar function $\eta = \eta (x)$ on $M$ such that
\be
r_{00}= \eta \alpha^{2}. \label{eq6.2}
\ee
Then (\ref{f2u}) and (\ref{formula23}) are reduced respectively to
\beq
&& \beta \zeta +b^{4}{\rm Ric}^{\alpha}+\left\{ \big[b^{k}\eta_{k}+(n-2)\eta^{2}\big]b ^{2}-u\right\}\alpha^{2}-(n-2)\eta^{2}\beta^{2} \nonumber\\
&& +\left\{(n-2)(b^{2}\eta_{0}-2\eta s_{0})+2(n-1)b^{2}\eta f_{0}\right\}\beta +(n-2)(b^{2}s_{0;0}-s_{0}^{2})\nonumber\\
&& +(n-1)b^{4}\left( {\rm Hess}_{F}f(y) +f_{0}^{2}\right)=0, \label{eq6.3}
\eeq
\beq
&& \left\{(n-2)s^{i}s_{i}+3(n-1)b^{2}\theta_{i}b^{i} - b^{2}(s^{k}_{\ ;k}+s^{k}_{ \ i}s^{i}_{\ k})\right\}\beta +(n-3)b^{2}\eta s_{0} \nonumber\\
&& +b^{2}b^{k}s_{0;k}-b^{4}s^{k}_{\ 0;k}+(n-1)b^{2}s_{k}s^{k}_{\ 0}-3(n-1)b^{4}\theta =0. \label{eq6.4}
\eeq
In this case,
\be
u= (n-2)s^{i}s_{i}+3(n-1)b^{2}\theta_{i}b^{i}-b^{2}(s^{k}_{ \ ;k}+s^{k}_{\ i}s^{i}_{\ k}). \label{eq6.5}
\ee

Then we have proved the following

\begin{thm} \label{mainth5}
Let $F=\frac{\alpha^{2}}{\beta}$ be a Kropina metric on an n-dimensional manifold $M$. Then F is  weakly weighted Einstein-Kropina metric  satisfying
\[
{\rm PRic}=(n-1)\left(\frac{3\theta}{F}+\sigma\right)F^{2}
\]
with respect to some volume form $dV=e^{-(n+1)f}dV_{BH}$ if and only if equations  (\ref{eq6.2}) and (\ref{eq6.3}), (\ref{eq6.4})  hold for some 1-form $\zeta$ and $u= u(x)$ and $\sigma$ are determined by (\ref{eq6.5}) and (\ref{formula24}) respectively.

\end{thm}

\vskip 3mm

\noindent
Xinyue Cheng \\
School of Mathematical Sciences \\
Chongqing Normal University \\
Chongqing  401331,  P. R. of China  \\
E-mail: chengxy@cqnu.edu.cn

\vskip 4mm

\noindent
Hong Cheng \\
School of Mathematical Sciences \\
Chongqing Normal University \\
Chongqing  401331,  P. R. of China  \\
E-mail: 1316011134@qq.com

\vskip 4mm

\noindent
Pengsheng Wu \\
School of Mathematical Sciences \\
Chongqing Normal University \\
Chongqing  401331,  P. R. of China  \\
E-mail: wu.pengsheng@qq.com

\end{document}